\documentclass[10pt]{book}
\usepackage{array}
\usepackage{latexsym}
\usepackage{amsmath}
\usepackage{amsfonts}
\usepackage{threeparttable}
\usepackage[dvips]{graphicx}
\usepackage{setspace}
\newcommand {\apgt} {\ {\raise-.5ex\hbox{$\buildrel>\over\sim$}}\ }
\newcommand {\aplt} {\ {\raise-.5ex\hbox{$\buildrel<\over\sim$}}\ }
\newtheorem{theorem}{\indent \sc Theorem}
\newtheorem{corollary}{\indent \sc Corollary}

\newtheorem{remark}{\indent \sc Remark}

\begin{document}


\renewcommand{\baselinestretch}{1.2}

\markright{ \hbox{\footnotesize\rm (2007): }\hfill }

\markboth{\hfill{\footnotesize\rm $K$-Step Maximum Likelihood Estimate} \hfill} {\hfill {\footnotesize\rm $K$-Step Maximum Likelihood Estimate}
\hfill}

\renewcommand{\thefootnote}{}
$\ $\par


\fontsize{10.95}{14pt plus.8pt minus .6pt}\selectfont \vspace{0.8pc} \centerline{\large\bf Convergence Rate of $K$-Step Maximum Likelihood
Estimate } \vspace{2pt} \centerline{\large\bf in Semiparametric Models} \vspace{.4cm} \centerline{Guang Cheng} \vspace{.4cm} \centerline{\it
Duke University} \vspace{.55cm} \fontsize{9}{11.5pt plus.8pt minus .6pt}\selectfont


\begin{quotation}
\noindent {\it Abstract:} We suggest an iterative approach to computing $K$-step maximum likelihood estimates (MLE) of the parametric components
in semiparametric models based on their profile likelihoods. The higher order convergence rate of $K$-step MLE mainly depends on the precision
of its initial estimate and the convergence rate of the nuisance functional parameter in the semiparametric model. Moreover, we can show that
the $K$-step MLE is as asymptotically efficient as the regular MLE after a finite number of iterative steps. Our theory is verified for several
specific semiparametric models. Simulation studies are also presented to support these theoretical results. \par

\vspace{9pt} \noindent {\it Key words and phrases:} $K$-Step Maximum Likelihood Estimate, Convergence Rate, Profile Likelihood, Semiparametric
Models.
\par
\end{quotation}\par


\fontsize{10.95}{14pt plus.8pt minus .6pt}\selectfont

\setcounter{chapter}{1}
\setcounter{equation}{0} 
\noindent {\bf 1. Introduction}

Let $X_{1},\ldots,X_{n}$ be independent and identically distributed random variables from a semiparametric model
$\mathbf{P}=\{P_{\theta,\eta}:\theta\in\Theta,\eta\in\mathcal{H}\}$, where $\theta$ is a $d-$dimensional parameter of interest and $\eta$ is an
infinite dimensional nuisance parameter. A well-known method of estimating the parameter $\theta$ in a semiparametric model is to solve $\theta$
from the below estimation equation:
\begin{eqnarray}
\sum_{i=1}^{n}\tilde{\ell}_{\theta,\tilde{\eta}_{n}}(X_{i})=0,\label{estequ}
\end{eqnarray}
where $\tilde{\eta}_{n}$ is some estimator for the nuisance parameter, and $\tilde{\ell}_{\theta,\eta}$ is the efficient score function for
$\theta$, whose definition will be introduced later. 
However, there are at least two concerns in solving (\ref{estequ}). Firstly, we may have multiple roots in which identifying the consistent
solution could be very challenging. Secondly, the above estimation approach requires an explicit form of the efficient score function, which in
general is implicitly defined as an orthogonal projection. Although we can estimate $\theta$ only by solving
$\sum\dot{\ell}_{\theta,\eta_{0}}(X_{i})=0$, where $\dot{\ell}_{\theta,\eta_{0}}$ is the regular score function for $\theta$ given the true
parameter $\eta_{0}$, in the semiparametric models of convex parametrization (page 305 in Bickel, Klaassen, Ritov and Wellner (1998)), many
other semiparametric models of interest do not possess such nice properties.

The above concerns can be addressed well by the profile likelihood based $K$-step maximum likelihood estimate proposed in this paper. Under
fairly general assumptions the $K$-step MLE is shown to posses higher order asymptotic efficiency than MLE of $\theta$ in semiparametric models.
Actually the motivation for constructing $k$-step estimator $\hat{\theta}_{n}^{(k)}$ comes from the Newton-Raphson algorithm for solving
(\ref{estequ}) with respect to $\theta$, starting at the initial guess $\hat{\theta}_{n}^{(0)}$. Thus, we can define $k$-step estimator
iteratively in the below form:
\begin{eqnarray}
\hat{\theta}_{n}^{(k)}=\hat{\theta}_{n}^{(k-1)}+\left(\mathbb{P}_{n}\tilde{\ell}_{\hat{\theta}_{n}^{(k-1)},
\tilde{\eta}_{n}}\tilde{\ell}^{T}_{\hat{\theta}_{n}^{(k-1)}, \tilde{\eta}_{n}}
\right)^{-1}\mathbb{P}_{n}\tilde{\ell}_{\hat{\theta}_{n}^{(k-1)},\tilde{\eta}_{n}}\label{1step}
\end{eqnarray}
for $k=1,2,\dots$, $\mathbb{P}_{n}f=\sum_{i=1}^{n}f(X_{i})/n$ and $\hat{\theta}_{n}^{(0)}$ is some preliminary estimator for $\theta$. In the
parametric models, $K$-step MLE is defined similarly but with the efficient score function replaced by the regular score function for $\theta$
in (\ref{1step}). Under some regularity conditions in parametric models, Jassen, Jureckova and Veraverbeke (1985) shows that
\begin{eqnarray}
\hat{\theta}_{n}^{(1)}-\hat{\theta}_{n}=O_{P}(n^{-1})\;\;\mbox{and}\;\;\;\hat{\theta}_{n}^{(2)}-\hat{\theta}_{n}=O_{P}(n^{-3/2}),\label{pararate}
\end{eqnarray}
where $\hat{\theta}_{n}$ is maximum likelihood estimate for $\theta$. The previous studies (Bickel, Klaassen, Ritov and Wellner (1998) and Van
der Vaart (1998)) about $K$-step MLE only focus on the semiparametric models with convex parametrization, in which the efficient score functions
can be estimated explicitly. Given certain no-bias conditions of the estimated efficient score functions, Van der Vaart (1998) shows that
$\hat{\theta}_{n}^{(1)}=\hat{\theta}_{n}+o_{P}(n^{-1/2})$. Moreover, $K$-step approach is also used in local (quasi) likelihood estimation for
the purpose of reducing computational cost, see Fan and Chen (1999), Fan, Chen and Zhou (2006) and Cai, Fan and Li (2000). However, as far as we
are aware, it appears that no systematic studies have been done on the construction of $K$-step semiparametric MLE and its higher order
asymptotic efficiency so far.

The efficient score function $\tilde{\ell}_{\theta,\eta}$ in (\ref{1step}) usually does not have an explicit form or cannot be estimated
explicitly as discussed above. Hence, we estimate $\mathbb{P}_{n}\tilde{\ell}_{\theta,\tilde{\eta}_{n}}$ and
$\mathbb{P}_{n}\tilde{\ell}_{\theta,\tilde{\eta}_{n}}\tilde{\ell}_{\theta,\tilde{\eta}_{n}}^{T}$ via numerical derivatives of the profile
likelihood. The profile likelihood $pl_{n}(\theta)$ is defined as $\sup_{\eta\in\mathcal{H}}lik_{n}(\theta,\eta)$, where $lik_{n}(\theta,\eta)$
is the full likelihood given $n$ observations. In practice, the profile likelihood may have an explicit form, e.g. the Cox model with right
censored data, or can be easily computed using procedures such as the fixed-point algorithm (as used in~Kosorok, Lee and Fine (2004), for
example) or the iterative convex minorant algorithm introduced in Groeneboom (1991) if $\eta$ is a monotone function. Hence we will assume
throughout this paper that evaluation of $pl_n(\theta)$ is computationally feasible. We shall consider the profile likelihood based $K$-step MLE
in the form:
\begin{eqnarray}
\hat{\theta}_{n}^{(k)}&=&\hat{\theta}_{n}^{(k-1)}+\left(\Pi_{n}(\hat{\theta}_{n}^{(k-1)},t_{n})\right)^{-1}\Gamma_{n}(\hat{\theta}_{n}^{(k-1)},s_{n})
\label{sche}
\end{eqnarray}
for $k=1,2,\ldots$ and reasonably accurate starting point $\hat{\theta}_{n}^{(0)}$. $\Gamma_{n}(\theta,s_{n})$ and $\Pi_{n}(\theta,t_{n})$ are
thus the discretized version of first and second derivative of the profile likelihood around $\theta$ with step size $s_{n}$ and $t_{n}$,
respectively. Their forms are given and justified in section 3. In section 2, we provide some necessary background about semiparametric models
and two primary assumptions needed in this paper. In section 3, we discuss the construction of the initial estimates and present the main result
of the paper about higher order convergence rate of $K$-step semiparametric MLE. In section 4, the proposed $K$-step approach is applied to
three semiparametric models. Section 5 contains some simulations results of the Cox regression model, and proofs are given in section~6.
\par

\vspace{0.2 in} \setcounter{chapter}{2}
\setcounter{equation}{0} 
\noindent {\bf 2. Background and Assumptions}

We assume the data $X_1,\ldots,X_n$ are i.i.d. throughout the paper. In what follows, we first briefly review the concept of the efficient score
function and define the convergence rate for the nuisance functional parameter. Next, we present two primary assumptions about second order
asymptotic expansions of log-profile likelihood and MLE.

\noindent {\bf 2.1 Preliminary}

The {\it score function} for $\theta$, $\dot{\ell}_{\theta,\eta}$, is defined as the partial derivative w.r.t. $\theta$ of the log-likelihood
given $\eta$ is fixed for a single observation. We denote the true values of $(\theta,\eta)$ as $(\theta_{0},\eta_{0})$. A score function for
$\eta_{0}$ is of the form
\begin{displaymath}
\frac{\partial}{\partial t}|_{t=0}\log p_{\theta_{0},\eta_{t}}(x)\equiv A_{\theta_{0},\eta_{0}}h(x),
\end{displaymath}
where $h$ is a ``direction'' by which $\eta_{t}\in\mathcal{H}$ approaches $\eta_{0}$, running through some index set H.
$A_{\theta,\eta}:H\mapsto L_{2}^{0}(P_{\theta,\eta})$ is the score operator for $\eta$. The {\it efficient score function} for $\theta$ is
defined as $\tilde{\ell}_{\theta,\eta}= \dot{\ell}_{\theta,\eta}-\Pi_{\theta,\eta}\dot{\ell}_{\theta,\eta}$, where
$\Pi_{\theta,\eta}\dot{\ell}_{\theta,\eta}$ minimizes the squared distance $P_{\theta,\eta}(\dot{\ell}_{\theta,\eta}-k)^{2}$ over all functions
$k$ in the closed linear space of the score functions for $\eta$ (the ``nuisance scores''). The inverse of the variance of
$\tilde{\ell}_{\theta,\eta}$ is the Cr\'{a}mer Rao bound for estimating $\theta$ in the presence of the infinite dimensional nuisance parameter
$\eta$, called efficient information matrix $\tilde{I}_{\theta,\eta}$. We also abbreviate $\tilde{\ell}_{\theta_{0},\eta_{0}}$ and
$\tilde{I}_{\theta_{0},\eta_{0}}$ with $\tilde{\ell}_{0}$ and $\tilde{I}_{0}$, respectively. An insightful review of efficient score functions
can be found in chapter~3 of Kosorok (2007).

The maximum likelihood estimate for $(\theta,\eta)$ can be expressed as $(\hat{\theta}_{n},\hat{\eta}_{n})$, where
$\hat{\eta}_n=\hat{\eta}_{\hat{\theta}_n}$ and $\hat{\eta}_{\theta}=argmax_{\eta\in\mathcal{H}} lik_{n}(\theta,\eta)$. The convergence rate for
$\eta$ is defined as the largest $r$ that satisfies
$\|\hat{\eta}_{\tilde{\theta}_{n}}-\eta_{0}\|=O_{P}(\|\tilde{\theta}_{n}-\theta_{0}\|+n^{-r})$, where $\|\cdot\|$ is a norm with definition
depending on context, i.e., for a Euclidean vector $u$, $\|u\|$ is the Euclidean norm, and for an element of the nuisance parameter space
$\eta\in{\cal H}$, $\|\eta\|$ is some chosen norm on ${\cal H}$. In regular semiparametric models, which we can define without loss of
generality to be models where the entropy integral converges, $r$ is always larger than $1/4$. We say the nuisance parameter has parametric rate
if $r=1/2$. For instance, the nuisance parameters of the three examples in Cheng and Kosorok (2006) achieve the parametric rate. More
specifically, the nuisance parameter in the Cox model, which is the cumulative hazard function, has the parametric rate under right censored
data. However, the convergence rate for the cumulative hazard becomes slower, i.e. $r=1/3$, under current status data.

\noindent {\bf 2.2 Assumptions}

The main result of this paper is based on the following second order asymptotic expansion of the profile likelihood, i.e. (\ref{lnplexp}). For
any random sequence $\tilde{\theta}_{n}=\theta_{0}+o_{P}(1)$, Cheng and Kosorok (2007) proves that
\begin{eqnarray}
log pl_{n}(\tilde{\theta}_{n})&=&\log pl_{n}(\theta_{0})+(\tilde{\theta}_{n}-\theta_{0})^{T} \sum_{i=1}^{n}\tilde{\ell}_{0}(X_{i})
\nonumber\\&&-\frac{n}{2}(\tilde{\theta}_{n}-\theta_{0})^{T}\tilde{I}_{0}(\tilde{\theta}_{n}-\theta_{0})
+O_P(g_{r}(\|\tilde{\theta}_{n}-\hat{\theta}_{n}\|)),\label{lnplexp}
\end{eqnarray}
where $g_{r}(w)\equiv(nw^{3}\vee n^{1-2r}w\vee n^{-2r+1/2})1\{1/4<r< 1/2\}+(nw^{3}\vee n^{-1/2})1\{r\geq 1/2\}$, under certain second order
no-bias conditions. Under similar conditions the maximum likelihood estimate is asymptotically normal, and has the asymptotic expansion:
\begin{eqnarray}
\sqrt{n}(\hat{\theta}_{n}-\theta_{0})&=& \frac{1}{\sqrt{n}}\sum_{i=1}^{n}\tilde{I}_{0}^{-1}\tilde{\ell}_{0}(X_{i}) +O_P(n^{-1/2}+n^{-2r+1/2}),
\label{mleexp}
\end{eqnarray}
where $\tilde{I}_{0}$ is assumed to be strictly positive definite. Expansions~(\ref{lnplexp}) and~(\ref{mleexp}) are essentially second order
versions of (1.4) and~(1.5), which justify using a semiparametric profile likelihood as an ordinary likelihood, in Murphy and Van der Vaart
(2000). Under second order conditions specified in section 2.3 of Cheng and Kosorok (2007), (\ref{lnplexp}) and (\ref{mleexp}) have been shown
to hold in several semiparametric models, e.g. Cox regression and partly linear model, in Cheng and Kosorok (2006) and Cheng and Kosorok (2007).
Therefore, we assume (\ref{lnplexp}) and (\ref{mleexp}) as two primary assumptions needed for the remainder of the paper.

\par

\vspace{0.2 in}

\setcounter{chapter}{3}
\setcounter{equation}{0} 
\noindent {\bf 3. Main Results}

We first present two general approaches to searching for the preliminary estimates. And then we discuss how to construct the estimates for
$\mathbb{P}_{n}\tilde{\ell}_{\theta,\tilde{\eta}_{n}}$ and
$\mathbb{P}_{n}\tilde{\ell}_{\theta,\tilde{\eta}_{n}}\tilde{\ell}_{\theta,\tilde{\eta}_{n}}^{T}$ in (\ref{sche}) based on the profile
likelihoods. Finally the convergence rate of $K$-step MLE is given. Such higher order convergence rate results are of interest particularly in
small- or moderate-sized samples. The conditions (\ref{lnplexp}) and (\ref{mleexp}) are assumed to hold in this section.

\noindent {\bf 3.1 Initial Estimate}

The start-up estimator is usually required to have reasonably good precision in the above $K$-step approach. In the parametric models,
$\hat{\theta}_{n}^{(0)}$ is required to be $\sqrt{n}$ consistent such that one- and two-step MLE can achieve the convergence rate as shown in
(\ref{pararate}). In our semiparametric model set-up, we need the initial estimate to be $n^{\psi}$ consistent for $0<\psi\leq 1/2$. The
$\sqrt{n}$ consistent estimate in parametric models can be obtained through M-estimation theorem, i.e. theorem 5.21 in Van der Vaart (1998), or
derived case by case in different examples. In the semiparametric models where the ad-hoc estimation methods for $\hat{\theta}_{n}^{(0)}$ are
unavailable, we provide two general search strategies for $\hat{\theta}_{n}^{(0)}$: one is through some MCMC sampling procedure, called the
profile sampler Lee, Kosorok and Fine (2005); another is through the deterministic or stochastic grid search over the profile likelihood
function.

The profile sampler is the MCMC sampling from the posterior of the profile likelihood, and was proposed for the purpose of obtaining frequentist
inference of $\theta$ Lee, Kosorok and Fine (2005). However, here we can use this convenient MCMC sampling procedure to yield
$\sqrt{n}$-consistent $\hat{\theta}_{n}^{(0)}$ and consistent estimate for $\tilde{I}_{0}$. Specifically speaking, under the conditions (1.4),
(1.5) in Murphy and Van der Vaart (2000) and mild conditions on the prior specified in theorem 1 of Lee, Kosorok and Fine (2005), Lee, Kosorok
and Fine (2005) shows that
\begin{eqnarray}
\tilde{E}_{\theta|\tilde{X}}(\theta)&=&\hat{\theta}_{n}+o_{P}(n^{-1/2}),\label{postmean}\\
\hat{I}_{n}(PS)&=&\tilde{I}_{0}+o_{P}(1),\label{postvar}
\end{eqnarray}
where $\tilde{E}_{\theta|\tilde{X}}(\theta)$ and $\hat{I}_{n}(PS)$ are the sample mean and the inverse of the sample variance of the profile
sampler, respectively.

Next, we provide an alternative grid search method to establish the $n^{\psi}$ consistent start-up estimator when the above profile sampling
procedure is unavailable or time consuming. When the dimension of $\theta$ is not large, we will conduct a deterministic search of objective
function $Q_{n}(\theta)$, which is defined as $(logpl_{n}(\theta)/n)$, at regularly spaced grid over the whole compact parameter space $\Theta$.
We summarize this idea in the below theorem~\ref{init}. Meanwhile, we need to assume the asymptotic uniqueness of $\hat{\theta}_{n}$:
\begin{eqnarray}
Q_{n}(\tilde{\theta}_{n})-Q_{n}(\hat{\theta}_{n})=o_{P}(1)\;\;\;\mbox{implies}\;\;\;\tilde{\theta}_{n}-\theta_{0}=o_{P}(1),\label{asyuni}
\end{eqnarray}
for any random sequence $\{\tilde{\theta}_{n}\}\in\Theta$.

\begin{theorem}\label{init} Let $\mathcal{D}_{n}$ be a set of points $\theta_{i}^{D}$ regularly spaced throughout $\Theta$ with
cardinality larger than $cn^{d\psi}$ for some $c>0$. Suppose that the parameter space $\Theta$ be a compact subset of $\mathbb{R}^{d}$ and
(\ref{asyuni}) holds, then we have for $0<\psi\leq 1/4$
\begin{eqnarray}
\theta_{n}^{D}-\theta_{0}=O_{P}(n^{-\psi}),
\end{eqnarray}
where $\theta_{n}^{D}=argmax_{\mathcal{D}_{n}}Q_{n}(\theta)$ .
\end{theorem}
However, if the dimension $d$ is very large, we prefer the outcome of a stochastic search whose search points are formed by the realizations of
an independent random variable $\bar{\theta}$ with strictly positive density around $\theta_{0}$.

\begin{corollary}\label{stosea}
Assume that $\bar{\theta}$ is independent of $Q_{n}(\theta)$ for all $\theta\in\Theta$ and admits a density having support $\Theta$ and bounded
away from zero in some neighborhood of $\theta_{0}$. Let $\mathcal{S}_{n}$ be a set of independent copies of $\bar{\theta}$ with cardinality
larger than $cn^{2\psi}$ for some $c>0$. Suppose that the parameter space $\Theta$ be a compact subset of $\mathbb{R}^{d}$ and (\ref{asyuni})
holds, then we have for $0<\psi\leq 1/4$
\begin{eqnarray}
\theta_{n}^{S}-\theta_{0}=O_{P}(n^{-\psi}),
\end{eqnarray}
where $\theta_{n}^{S}=argmax_{\mathcal{S}_{n}}Q_{n}(\theta)$ .
\end{corollary}




\noindent {\bf 3.2 $K$-step MLE}

Before proceeding to give the convergence rate of $K$-step MLE, we first specify the forms of $\Gamma_{n}(\theta,s_{n})$ and $\Pi(\theta,t_{n})$
in (\ref{sche}). The intuitive idea behind the constructions of the estimators for $\mathbb{P}_{n}\tilde{\ell}_{\theta,\tilde{\eta}_{n}}$ and
$\mathbb{P}_{n}\tilde{\ell}_{\theta,\tilde{\eta}_{n}}\tilde{\ell}_{\theta,\tilde{\eta}_{n}}^{T}$ is to use $\hat{\eta}_{\theta}$ as
$\tilde{\eta}_{n}$ when making inferences about $\theta$.

Specifically speaking, the $i$th component of $\Gamma_{n}(\theta,s_{n})$ is constructed in the form:
\begin{eqnarray}
[\Gamma_{n}(\theta,s_{n})]_{i}&=&\mathbb{P}_{n}\left\{
\frac{\log lik(\theta+s_{n}v_{i},\hat{\eta}_{\theta+s_{n}v_{i}})-\log lik(\theta,\hat{\eta}_{\theta})}{s_{n}}\right\}\nonumber\\
&=&\frac{\log pl_{n}(\theta+s_{n}v_{i})-\log pl_{n}(\theta)}{ns_{n}},\label{estesco}
\end{eqnarray}
where step size $s_{n}\overset{P}{\rightarrow} 0$ and $v_{i}$ denotes the ith unit vector in $\mathbb{R}^{d}$.  Following similar logic, we can
define the $(i,j)$-th component of $\Pi_{n}(\theta,t_{n})$ as:
\begin{eqnarray}
[\Pi_{n}(\theta,t_{n})]_{i,j}=&-&\frac{\log pl_{n}(\theta+v_{i}t_{n}+v_{j}t_{n})+\log pl_{n}(\theta)}{nt_{n}^{2}}\nonumber\\&+&\frac{\log
pl_{n}(\theta+v_{i}t_{n})+\log pl_{n}(\theta+v_{j}t_{n})}{nt_{n}^{2}},\label{estei}
\end{eqnarray}
where step size $t_{n}\overset{P}{\rightarrow} 0$. (\ref{estei}) is also called observed profile information in Murphy and Van der Vaart (1999).
The lemma~1 in the appendix justifies the use of (\ref{estesco}) and (\ref{estei}) as consistent estimates of $\mathbb{P}_{n}\tilde{\ell}_{0}$
and $\tilde{I}_{0}$, respectively.

The convergence rate of $K$-step MLE  is certainly determined by the order of the step sizes in numerical differentiations
$\Gamma_{n}(\cdot,s_{n})$ and $\Pi_{n}(\cdot,t_{n})$ as shown in the above. However, we are mostly interested in the fastest convergence rate
$K$-step MLE can attain. Hence, we assume using the optimal step sizes $(s_{n}^{\ast},t_{n}^{\ast})$, under which the fastest convergence rate
of $\hat{\theta}_{n}^{(k)}$ is achieved, in the below theorem~\ref{thm1}. As the theoretical basis for using $K$-step approach in practice, the
below theorem~\ref{thm1} first presents the convergence rate for the fully iterative estimate $\hat{\theta}_{n}^{(\infty)}$, called optimal rate
of $K$-step MLE, and then gives the number of iterations needed in (\ref{sche}) for $\hat{\theta}_{n}^{(k)}$ to attain the above optimal rate.
Note that the asymptotic efficiency of $\hat{\theta}_{n}^{(k)}$ has continuously improved through the whole iterative procedure until it reaches
the optimal bound based on the proof of theorem~\ref{thm1}.

\begin{theorem}\label{thm1}
Assume that $\hat{\theta}_{n}^{(k)}$ is defined as (\ref{sche}) and $\hat{\theta}_{n}^{(0)}$ is $n^{\psi}$-consistent for $0<\psi\leq 1/2$, we
have
\begin{eqnarray}
\hat{\theta}_{n}^{(\infty)}-\hat{\theta}_{n}=O_{P}(n^{-3/4}\vee n^{-r-1/4}).\label{optrate}
\end{eqnarray}
Moreover, the above optimal rate can be achieved after $N$ $(M)$ iterations starting from $\hat{\theta}_{n}^{(0)}$ in (\ref{sche}) for $r\geq
1/2$ $(1/4<r<1/2)$:
\begin{eqnarray}
\hat{\theta}_{n}^{(N)}-\hat{\theta}_{n}&=&O_{P}(n^{-3/4}),\label{1-mle}\\
\hat{\theta}_{n}^{(M)}-\hat{\theta}_{n}&=&O_{P}(n^{-r-1/4}),\label{2-mle}
\end{eqnarray}
where $N=int[\log 2\psi/\log(2/3)]+1$, $M=int[\log(\psi/r)/\log(2/3)]+int[\log(4r/(4r-1))/\log(2)-1]+1$ and $int[x]$ indicates the smallest
nonnegative integer larger than or equal to $x$.
\end{theorem}


From theorem~\ref{init} and~\ref{thm1}, it is not surprising to find that there exists a tradeoff between the number of search grids and the
number of iterations. Combining (\ref{1-mle}) and (\ref{2-mle}) with (\ref{mleexp}), we have the following asymptotic expansion of $K$-step MLE:
\begin{eqnarray*}
\sqrt{n}(\hat{\theta}_{n}^{(N)}-\theta_{0})&=& \frac{1}{\sqrt{n}}\sum_{i=1}^{n}\tilde{I}_{0}^{-1}\tilde{\ell}_{0}(X_{i})+O_P(n^{-1/4}),\\
\sqrt{n}(\hat{\theta}_{n}^{(M)}-\theta_{0})&=& \frac{1}{\sqrt{n}}\sum_{i=1}^{n}\tilde{I}_{0}^{-1}\tilde{\ell}_{0}(X_{i})+O_P(n^{-r+1/4}).
\end{eqnarray*}

Thus we can construct the $(1-\alpha)$-th two sided asymptotically correct confidence interval for $\theta$ based on $K$-step MLE, i.e.
$(\hat{\theta}_{n}^{(k)}-z_{1-\alpha/2}/\sqrt{n\tilde{I}},\hat{\theta}_{n}^{(k)}+z_{1-\alpha/2}/\sqrt{n\tilde{I}})$, where $k=M$ or $N$,
$z_{\alpha}$ is the standard normal $\alpha$-th quantile, and $\tilde{I}$ is a consistent estimator of $\tilde{I}_{0}$.


\begin{remark}
Recall that in the parametric models, Jassen, Jureckova and Veraverbeke (1985) shows that
$\hat{\theta}_{n}^{(1)}-\hat{\theta}_{n}=O_{P}(n^{-1})$ and $\hat{\theta}_{n}^{(2)}-\hat{\theta}_{n}=O_{P}(n^{-3/2})$. However, the optimal rate
for the $K$-step MLE is slower even in the semiparametric models with parametric convergence rate. Such efficiency loss can be partially
explained by the less smoothness of the profile likelihood in semiparametric models. In other words, the corresponding estimators for the score
function and information matrix in $K$-step parametric MLE usually have bias of smaller order.
\end{remark}

\par

\setcounter{chapter}{4}
\setcounter{equation}{0} 
\noindent {\bf 4. Examples}

In this section, the above $K$-step estimation approach is illustrated with three semiparametric models of different convergence rates. Under
the model assumptions specified in section 5 of Cheng and Kosorok (2007), Cheng and Kosorok (2007) shows that (\ref{lnplexp}) and (\ref{mleexp})
hold in all the examples. Hence, we only briefly review the model set-up here, and then discuss the choices of the initial estimates. Finally,
we apply the theorem~\ref{thm1} to figure out the least number of iterations in $K$-step MLE needed to achieve the full efficiency.

\noindent {\bf 4.1 Cox regression with right censored data}\label{coxrt}

In the Cox regression model, the hazard function of the survival time $T$ of a subject with covariate $Z$ is expressed as:
\begin{eqnarray}
\lambda(t|z)\equiv\lim_{\Delta\rightarrow 0}\frac{1}{\Delta}Pr(t\leq T<t+\Delta|T\geq t,Z=z)=\lambda(t)\exp(\theta z),\label{eg1den}
\end{eqnarray}
where $\lambda$ is an unspecified baseline hazard function and $\theta$ is a vector including the regression parameters Cox (1972). Under right
censoring data, we only know that the event time $T$ has occurred either before the censoring time $C$, or after the censoring time $C$. More
precisely, the data observed is $X=(Y,\delta,Z)$, where $Y=T\wedge C$, $\delta=I\{T\leq C\}$, and $Z\in\mathbb{Z}\subset\mathbb{R}$ is a
regression covariate. In the Cox regression model, we are usually interested in the regression parameter $\theta$ while treating the cumulative
hazard function $\eta$ as the nuisance parameter. Thus we express the likelihood for $(\theta,\eta)$ in the below form:
\begin{eqnarray}
lik(\theta,\eta)=\left(e^{\theta z}\eta\{y\}e^{-e^{\theta z}\eta(y)}\right)^{\delta}\left(e^{-e^{\theta z}\eta(y)}\right)^{1-\delta},
\end{eqnarray}
by replacing hazard function $\lambda(y)$ by the point mass $\eta\{y\}$. By the special construction of the Cox model, we have an explicit form
of the log-profile likelihood:
\begin{eqnarray}
\log pl_{n}(\theta)=\sum_{i=1}^{l}(\theta z_{[i]}-\log\sum_{j\in R_{i}}e^{\theta z_{j}}),\label{pleg1}
\end{eqnarray}
where $R_{i}=\{j:Y_{j}\geq t_{i}\}$, $t_{i}$ is the observed value of the $i$-th ordered event time and $z_{[i]}$ is the covariate corresponding
to $t_{i}$. The convergence rate of the estimated nuisance parameter is established in theorem 3.1 of Murphy and Van der Vaart (1999):
\begin{eqnarray}
\|\hat{\eta}_{\tilde{\theta}_{n}} -\eta_{0}\|_{\infty}=O_P(n^{-\frac{1}{2}}+ \|\tilde{\theta}_{n}-\theta_{0}\|),\label{eg1crate}
\end{eqnarray}
where $\|\cdot\|_{\infty}$ denotes the uniform norm.

In this model, the profile sampler is generated very fast because of the explicit form for the profile likelihood. Hence, we use it to yield the
root-$n$ consistent start-up estimator. By theorem~\ref{thm1}, we can conclude that $\hat{\theta}_{n}^{(1)}-\hat{\theta}_{n}=O_{P}(n^{-3/4})$,
where $\hat{\theta}_{n}^{(1)}$ is constructed according to (\ref{sche}).

\noindent {\bf 4.2 Cox regression for current status data}\label{coxcs}

Current status data arises when each subject is observed at a single examination time, $Y$, to determine if an event has occurred. The event
time, $T$, cannot be known exactly. Then the observed data are $n$ i.i.d. realizations of $X=(Y, \delta, Z)\in R^{+}\times \lbrace 0,1
\rbrace\times R$, where $\delta=I\{T \leq Y\}$ and $Z$ is a vector of covariates. It is not difficult to derive the log-likelihood:
\begin{eqnarray}
\log lik_{n}(\theta,\eta)=\sum_{i=1}^{n}\delta_{i}\log[1-\exp(-\eta(Y_{i})\exp(\theta Z_{i}))]-(1-\delta_{i})\exp(\theta
Z_{i})\eta(Y_{i}).\label{eg1lik}
\end{eqnarray}
Moreover, using entropy methods, Murphy and Van der Vaart (1999) extends earlier results of Huang (1996), show that
\begin{eqnarray}
\|\hat{\eta}_{\tilde{\theta}_{n}}-\eta_{0}\|_{L_{2}}=O_{P}(\|\tilde{\theta}_{n}-\theta_{0}\|+n^{-1/3}), \label{eg1rate}
\end{eqnarray}
where $\|\cdot\|_{L_{2}}$ is the $L_{2}$ norm w.r.t. the distribution of $Y$.

In the Cox regression with current status data, the iterative convex minorant algorithm Huang (1996) is implemented to yield the profile
likelihood. The MCMC sampling procedure thus becomes more time consuming because of such iterative computation mechanism. Hence, we prefer using
grid search approach to obtain $n^{1/4}$-consistent preliminary estimate. We know that three step MLE attains the optimal rate, i.e.
$\hat{\theta}_{n}^{(3)}-\hat{\theta}_{n}=O_{P}(n^{-7/12})$, based on theorem~\ref{thm1}.

\noindent {\bf 4.3 The partly linear model}

In this model, a continuous outcome $Y$, conditional on the covariates $(W,Z)\in \mathbb{R}^{d}\times\mathbb{R}$, is modeled as:
\begin{eqnarray}
Y=\theta^{T}W+k(Z)+\xi,\label{eg2lik}
\end{eqnarray}
where $k$ is an unknown smooth function, and $\xi\sim N(0,1)$. The functional nuisance parameter $k$ is assumed to belong to
$\mathcal{O}_{2}\equiv\{f:J_{2}(f)+\|f\|_{\infty}< M,\;\mbox{for a known}\;M<\infty\}$, where $J_{2}(f)$ is the second order Sobolev norm of
$f$. However, the response $Y$ is not observed directly, but only its current status is observed at a random censoring time $C\in\mathbb{R}$. In
other words, we observe $X=(C,\Delta,W,Z)$, where $\Delta=1_{\{Y\leq C\}}$. Additionally $(Y,C)$ is assumed to be independent given $(W,Z)$.
Under the model~(\ref{eg2lik}), the log-likelihood for a single observation at $X=x\equiv(c,\delta,w,z)$ can be shown to have the form:
\begin{eqnarray}
loglik_{\theta,k}(x)=\delta\log\left\{\Phi\left(c-\theta w-k(z)\right)\right\}+(1-\delta)\log\left\{1-\Phi\left(c-\theta
w-k(z)\right)\right\},\label{eg2p}
\end{eqnarray}
where $\Phi$ is the standard normal distribution. In lemma 4 of Cheng and Kosorok (2007), we have shown
\begin{eqnarray}
\left\|\hat{k}_{\tilde{\theta}_{n}}-k_{0}\right\|_{L_{2}}&=&O_{P}(n^{-2/5}+\|\tilde{\theta}_{n}-\theta_{0}\|).\label{eg3ratre}
\end{eqnarray}
The rate $r=2/5$ is clearly faster than the cubic rate but slower than the parametric rate. Depending on the dimension of $\theta$, we can
choose the deterministic or random search for the starting estimate. Similarly, we can show that four iterations are needed to achieve the
optimal rate, i.e. $\hat{\theta}_{n}^{(4)}-\hat{\theta}_{n}=O_{P}(n^{-13/20})$, if $\hat{\theta}_{n}^{(0)}$ is $n^{1/4}$-consistent.

\noindent {\bf 5. Simulations}

It is of interest to see, at a finite sample, how good the $K$-step MLE is in comparison with the regular MLE. Hence, we conducted simulations
in the Cox regression model with right censored data and with current status data in this section. The simulation results presented in the table
1 and 2 agree with our theoretical results given in subsection~\ref{coxrt} and \ref{coxcs}.

We first run the simulations for various sample sizes in the Cox model with right censored data. As indicated in subsection~\ref{coxrt}, we can
construct $\hat{\theta}_{n}^{(1)}$ in the form of (\ref{sche}) with $(s_{n}^{\ast},t_{n}^{\ast})$ set to be proportional to
$(n^{-3/4},n^{-1/2})$ according to the proof of theorem~\ref{thm1}. The profile sampler is generated under a Lebesgue prior. For each sample
size, 500 datasets were analyzed. The event times were generated from (\ref{eg1den}) with one covariate $Z\sim U[0,1]$. The regression
coefficient is $\theta=1$ and $\eta(t)=\exp(t)-1$. The censoring time $C\sim U[0,t_{n}]$, where $t_{n}$ was chosen such that the average
effective sample size over 500 samples is approximately $0.9n$. For each dataset, Markov chains of length 5,000 with a burn-in period of 1,000
were generated using the Metropolis algorithm. The jumping density for the coefficient was normal with current iteration and variance tuned to
yield an acceptance rate of $20\%-40\%$. In the Cox regression with current status data, we first use the deterministic search over $[-5,5]$ for
the $n^{1/4}$ consistent $\hat{\theta}_{n}^{(0)}$. The three step MLE is iteratively generated according to (\ref{sche}), in which the order of
$(s_{n}^{\ast},t_{n}^{\ast})$ at each step is specified in the proof of theorem~\ref{thm1}.

In the appendix, the table 1 (2) summarizes the results from the simulations of Cox regression with right censored data (current status data)
giving the average across 500 samples of $K$-step MLE and the maximum likelihood estimate (MLE). According to theorem~\ref{thm1},
$n^{3/4}|\hat{\theta}_{n}-\hat{\theta}_{n}^{(1)}|$ ($n^{7/12}|\hat{\theta}_{n}-\hat{\theta}_{n}^{(3)}|$) in Cox model with right censored data
(current status data) is bounded in probability. And the realizations of these terms summarized in table 1 and 2 clearly illustrate their
boundedness. For each sample size, we can clearly observe that $K$-step MLE approaches to $\hat{\theta}_{n}$ after every iteration. Hence we can
conclude that the numerical evidence in this section supports our theoretical results.



\vspace{0.2 in} \noindent {\large\bf Acknowledgment} The author thank Dr. Michael Kosorok for several insightful discussions.
\par


\vspace{0.2 in}

\setcounter{chapter}{6}
\setcounter{equation}{0} 
\noindent {\bf 6. Appendix}

In the below lemma~1, we first provide a key technical tool for deriving higher order convergence rate of $K$-step MLE. The symbol $R_{n}\asymp
q_{n}$ means that some random quantity $R_{n}=O_{P}(q_{n})$ and $R_{n}^{-1}=O_{P}(q_{n}^{-1})$, where $q_{n}\rightarrow 0$.\\

{\it Lemma 1.} Assume the conditions (\ref{lnplexp}) and (\ref{mleexp}) and $\hat{\theta}_{n}^{(0)}$ is a $n^{\psi}$-consistent estimate for
$0<\psi\leq 1/2$, then we have
\begin{eqnarray}
\Gamma_{n}(\hat{\theta}_{n}^{(0)},s_{n})&=&\mathbb{P}_{n}\tilde{\ell}_{0}+O_{P}\left(n^{-\psi}\vee |s_{n}|\vee\frac{g_{r}(n^{-\psi}\vee|s_{n}|)}
{n|s_{n}|}\right),\label{ga1}\\
\Gamma_{n}(\hat{\theta}_{n}^{(0)}+U_{n},s_{n})&=&\Gamma_{n}(\hat{\theta}_{n}^{(0)},s_{n})-\tilde{I}_{0}U_{n}\nonumber\\&&+
O_{P}\left(|s_{n}|\vee\frac{g_{r}(n^{-\psi}\vee|s_{n}|\vee\|U_{n}\|)}{n|s_{n}|}\right),\label{lem2}\\
\Pi_{n}(\tilde{\theta}_{n},t_{n})&=&\tilde{I}_{0}+O_{P}\left(\frac{g_{r}(r_{n}\vee |t_{n}|)}{nt_{n}^{2}}\right),\label{estrel}
\end{eqnarray}
where $U_{n}=O_{P}(n^{-s})$ for some $s>0$, $(\tilde{\theta}_{n}-\hat{\theta}_{n})=O_P(r_{n})$ and $g_{r}(w)\equiv(nw^{3}\vee n^{1-2r}w\vee
n^{-2r+1/2})1\{1/4<r< 1/2\}+(nw^{3}\vee n^{-1/2})1\{r\geq 1/2\}$.

{\it Proof of lemma~1:} (\ref{lnplexp}) implies that
\begin{eqnarray*}
log pl_{n}(\hat{\theta}_{n}^{(0)}+V_{n}+s_{n}v_{i})&=&\log pl_{n}(\theta_{0})+(\hat{\theta}^{(0)}_{n}+V_{n}+s_{n}v_{i}-\theta_{0})^{T}
\sum_{i=1}^{n}\tilde{\ell}_{0}(X_{i})\\
&&-\frac{n}{2}(\hat{\theta}_{n}^{(0)}+V_{n}+s_{n}v_{i}-\theta_{0})^{T}\tilde{I}_{0}(\hat{\theta}_{n}^{(0)}+V_{n}+s_{n}v_{i}-\theta_{0})\\
&&+O_P(g_{r}(n^{-\psi}\vee|s_{n}|\vee\|V_{n}\|)),\\
log pl_{n}(\hat{\theta}_{n}^{(0)}+V_{n})&=&\log pl_{n}(\theta_{0})+(\hat{\theta}^{(0)}_{n}+V_{n}-\theta_{0})^{T}
\sum_{i=1}^{n}\tilde{\ell}_{0}(X_{i})\\
&&-\frac{n}{2}(\hat{\theta}_{n}^{(0)}+V_{n}-\theta_{0})^{T}\tilde{I}_{0}(\hat{\theta}_{n}^{(0)}+V_{n}-\theta_{0})\\&&+O_P(g_{r}(n^{-\psi}\vee\|V_{n}\|)),
\end{eqnarray*}
for any random vector $V_{n}=o_{P}(1)$ and $s_{n}\overset{P}{\rightarrow}0$. Combining the above two expansions and (\ref{estesco}), we have
\begin{eqnarray*}
\Gamma_{n}(\hat{\theta}_{n}^{(0)}+V_{n},s_{n})=\mathbb{P}_{n}\tilde{\ell}_{0}-\tilde{I}_{0}(\hat{\theta}_{n}^{(0)}-\theta_{0})-\tilde{I}_{0}V_{n}
+O_{P}\left(|s_{n}|\vee\frac{g_{r}(n^{-\psi}\vee|s_{n}|\vee\|V_{n}\|)}{n|s_{n}|}\right).
\end{eqnarray*}
By replacing $V_{n}=0$ and $V_{n}=U_{n}$ in the above equation, we have proved (\ref{ga1}) and (\ref{lem2}), respectively. Taking into account
(\ref{lnplexp}) and (\ref{mleexp}), we can prove the below second order asymptotic expansion of the profile likelihood around
$\hat{\theta}_{n}$:
\begin{eqnarray}
\log pl_{n}(\tilde{\theta}_{n})=\log pl_{n}(\hat{\theta}_{n})-\frac{1}{2}n(\tilde{\theta}_{n}-
\hat{\theta}_{n})^{T}\tilde{I}_{0}(\tilde{\theta}_{n}- \hat{\theta}_{n})+O_{P}(g_{r}(\|\tilde{\theta}_{n}-\hat{\theta}_{n}\|))
\label{lnplexphat}
\end{eqnarray}
for any sequence $\tilde{\theta}_{n}=\hat{\theta}_{n}+o_{P}(1)$. Following similar analysis in the above, (\ref{estei}) and (\ref{lnplexphat})
yield (\ref{estrel}). This completes the whole proof. $\Box$

{\it Proof of theorem~\ref{init}:} (\ref{lnplexp}) implies that for $\tilde{\theta}_{n}-\theta_{0}=o_{P}(1)$
\begin{eqnarray}
Q_{n}(\tilde{\theta}_{n})&=&Q_{n}(\theta_{0})+(\tilde{\theta}_{n}-\theta_{0})^{T}\mathbb{P}_{n}\tilde{\ell}_{0}-\frac{1}{2}(\tilde{\theta}_{n}
-\theta_{0})^{T}\tilde{I}_{0}(\tilde{\theta}_{n}-\theta_{0})+\Delta_{n},\label{inter1}
\end{eqnarray}
where $\Delta_{n}=O_{P}(g_{r}(\|\tilde{\theta}_{n}-\hat{\theta}_{n}\|))/n$. We then show that
$P(\|\theta_{n}^{D}-\theta_{0}\|>Cn^{-\psi})\rightarrow 0$ by the below set of inequalities for some $C>0$. Set
$\mathcal{N}_{n}=\{\theta:\|\theta-\theta_{0}\|\leq Cn^{-\psi}\}$ and $\mathcal{N}_{n}^{c}$ denotes its complement. Note that
$\mathcal{D}_{n}\cap\mathcal{N}_{n}\neq\emptyset$ for $C$ large enough and $\mathcal{D}_{n}\cap\mathcal{N}_{n}^{c}\neq\emptyset$ for $n$ large
enough.
\begin{eqnarray*}
P(\theta_{n}^{D}\in\mathcal{N}_{n}^{c}) &\leq& P\left(\max_{\mathcal{D}_{n}\cap\mathcal{N}_{n}}Q_{n}(\theta)
\leq\max_{\mathcal{D}_{n}\cap\mathcal{N}_{n}^{c}}Q_{n}(\theta)\right)\\
&\leq& P\left(\max_{\mathcal{D}_{n}\cap\mathcal{N}_{n}}Q_{n}(\theta)<
Q_{n}(\theta_{0})-C_{1}n^{-2\psi}\right)\\
&&+ P\left(\left\{\max_{\mathcal{D}_{n}\cap\mathcal{N}_{n}}Q_{n}(\theta)
\leq\max_{\mathcal{D}_{n}\cap\mathcal{N}_{n}^{c}}Q_{n}(\theta)\right\}\cap\right.\\
&&\mbox{\hspace{2in}}\left.\left\{\max_{\mathcal{D}_{n}\cap\mathcal{N}_{n}}Q_{n}(\theta)\geq
Q_{n}(\theta_{0})-C_{1}n^{-2\psi}\right\}\right)\\
&\leq&P\left(\max_{\mathcal{D}_{n}\cap\mathcal{N}_{n}}\sqrt{n}(Q_{n}(\theta)-Q_{n}(\theta_{0}))<-C_{1}n^{1/2-2\psi}\right)\\
&&+P\left(\max_{\mathcal{D}_{n}\cap\mathcal{N}_{n}^{c}}\sqrt{n}(Q_{n}(\theta)-Q_{n}(\theta_{0}))\geq-C_{1}n^{1/2-2\psi}\right),
\end{eqnarray*}
where $C_{1}$ is some positive number. The first inequality in the above follows from the definition of $\theta_{n}^{D}$. Based on
(\ref{inter1}) we have
\begin{eqnarray*}
&&P\left(\max_{\mathcal{D}_{n}\cap\mathcal{N}_{n}}\sqrt{n}(Q_{n}(\theta)-Q_{n}(\theta_{0}))<-C_{1}n^{1/2-2\psi}\right)\nonumber\\
&&=P\left(\max_{\mathcal{D}_{n}\cap\mathcal{N}_{n}}\left(\sqrt{n}(\theta-\theta_{0})^{T}\mathbb{P}_{n}\tilde{\ell}_{0}
-\frac{\sqrt{n}}{2}(\theta-\theta_{0})^{T}\tilde{I}_{0}(\theta-\theta_{0})+\sqrt{n}\Delta_{n}\right)
<-C_{1}n^{1/2-2\psi}\right)\nonumber\\
&&\leq P\left(\max_{\mathcal{D}_{n}\cap\mathcal{N}_{n}}(\theta-\theta_{0})(-\sqrt{n}\mathbb{P}_{n}\tilde{\ell}_{0})+\max_{\mathcal{D}_{n}
\cap\mathcal{N}_{n}}((\sqrt{n}/2)(\theta-\theta_{0})^{T}\tilde{I}_{0}(\theta-\theta_{0}))+\right.\\
&&\mbox{\hspace{3.5in}}\left.\max_{\mathcal{D}_{n}
\cap\mathcal{N}_{n}}(\sqrt{n}\Delta_{n})>C_{1}n^{1/2-2\psi}\right)\nonumber\\
&&\leq P\left(\sqrt{n}\mathbb{P}_{n}\tilde{\ell}_{0}\apgt (C_{1}-\delta C^{2}/2)n^{1/2-2\psi}+O_{P}(n^{1/2-3\psi})\right),
\end{eqnarray*}
where $\delta$ is the largest eigenvalue for $\tilde{I}_{0}$. The last inequality in the above follows from the compactness of $\Theta$. Let
$\theta_{n}^{\ast}=argmax_{\mathcal{D}_{n}\cap\mathcal{N}_{n}^{c}}\sqrt{n}Q_{n}(\theta)$. (\ref{asyuni}) implies that
$\theta_{n}^{\ast}-\theta_{0}=o_{P}(1)$ since $Q_{n}(\hat{\theta}_{n})-Q_{n}(\theta_{0})=o_{P}(1)$. Thus, by (\ref{inter1}), we have
\begin{eqnarray*}
&&P\left(\max_{\mathcal{D}_{n}\cap\mathcal{N}_{n}^{c}}\sqrt{n}(Q_{n}(\theta)-Q_{n}(\theta_{0}))\geq-C_{1}n^{1/2-2\psi}\right)\nonumber\\
&&=P\left(\sqrt{n}(\theta_{n}^{\ast}-\theta_{0})\mathbb{P}_{n}\tilde{\ell}_{0}-\frac{\sqrt{n}}{2}(\theta_{n}^{\ast}-\theta_{0})^{T}\tilde{I}_{0}
(\theta_{n}^{\ast}-\theta_{0})+O_{P}(g_{r}(\|\theta_{n}^{\ast}-\hat{\theta}_{n}\|)/\sqrt{n})\right.\\
&&\mbox{\hspace{4in}}\left.\geq -C_{1}n^{1/2-2\psi}\right)\nonumber\\
&&\leq P\left(\sqrt{n}\mathbb{P}_{n}\tilde{\ell}_{0}\apgt(\delta K_{1}^{2}/2-C_{1})n^{1/2-2\psi}+O_{P}(n^{1/2-3\psi})\right)
\end{eqnarray*}
Note that $\theta_{n}^{\ast}$ belongs to the regularly spaced grid set $\mathcal{D}_{n}$ and $\|\theta_{n}^{\ast}-\theta_{0}\|>Cn^{-\psi}$.
Therefore, we can conclude that $\theta_{n}^{\ast}$ should be the closest grid point to $\theta_{0}$ but not in $\mathcal{N}_{n}$, i.e.
$K_{1}n^{-\psi}\leq\|\theta_{n}^{\ast}-\theta_{0}\|\leq K_{2}n^{-\psi}$, where $C<K_{1}\leq K_{2}\leq 2C$ for large $C$, from (\ref{inter1}) and
the construction of $\mathcal{D}_{n}$. Without loss of generality, we assume $\theta_{n}^{\ast}>\theta_{0}$. Thus the last inequality in the
above follows. Note that $\sqrt{n}\mathbb{P}_{n}\tilde{\ell}_{0}=O_{P}(1)$ and $\psi\leq 1/4$. By choosing sufficiently large $C$ and $C_{1}$,
meanwhile keeping the inequality $\delta C^{2}/2<C_{1}<\delta K_{1}^{2}/2$ hold, we can $P(\theta_{n}^{D}\in\mathcal{N}_{n}^{c})\rightarrow 0$
based on the above inequalities. $\Box$

{\it Proof of corollary~\ref{stosea}:} The proof is similar to that of theorem~\ref{init}. We still need to show that
$P(\|\theta_{n}^{S}-\theta_{0}\|>Cn^{-\psi})\rightarrow 0$ for some $C>0$. Similarly, we have
\begin{eqnarray*}
P(\theta_{n}^{S}\in\mathcal{N}_{n}^{c})&\leq&E\left[P\left(\max_{\mathcal{S}_{n}\cap\mathcal{N}_{n}}Q_{n}(\theta)
\leq\max_{\mathcal{S}_{n}\cap\mathcal{N}_{n}^{c}}Q_{n}(\theta)|\mathcal{S}_{n}\right)\right]\\
&\leq&E\left[P\left(\max_{\mathcal{S}_{n}\cap\mathcal{N}_{n}}\sqrt{n}(Q_{n}(\theta)-Q_{n}(\theta_{0}))<-C_{1}n^{1/2-2\psi}|\mathcal{S}_{n}\right)\right]\\
&&+P\left(\sup_{\mathcal{N}_{n}^{c}}\sqrt{n}(Q_{n}(\theta)-Q_{n}(\theta_{0}))\geq-C_{1}n^{1/2-2\psi}\right)\\
&\leq& P\left(\sqrt{n}\mathbb{P}_{n}\tilde{\ell}_{0}\apgt (C_{1}/2)n^{1/2-2\psi}+O_{P}(n^{1/2-3\psi})\right)\\&&+
P\left(\sup_{\mathcal{N}_{n}^{c}}\sqrt{n}(Q_{n}(\theta)-Q_{n}(\theta_{0}))\geq-C_{1}n^{1/2-2\psi}\right)\\
&&+E\left[P\left(\min_{\mathcal{S}_{n}\cap
\mathcal{N}_{n}}((\sqrt{n}/2)(\theta-\theta_{0})^{T}\tilde{I}_{0}(\theta-\theta_{0}))>(C_{1}/2)n^{1/2-2\psi}|\mathcal{S}_{n}\right)\right]
\end{eqnarray*}
The first two quantity in the last inequality approaches to zero by choosing proper $C_{1}$ and $C$ according to similar analysis in the proof
of theorem~\ref{init}. We next analyze the last quantity.
\begin{eqnarray*}
&&E\left[P\left(\min_{\mathcal{S}_{n}\cap
\mathcal{N}_{n}}((\sqrt{n}/2)(\theta-\theta_{0})^{T}\tilde{I}_{0}(\theta-\theta_{0}))>(C_{1}/2)n^{1/2-2\psi}|\mathcal{S}_{n}\right)\right]\\
&&\leq E\left[P\left(\min_{\mathcal{S}_{n}}((\theta-\theta_{0})^{T}\tilde{I}_{0}(\theta-\theta_{0}))>C_{1}n^{-2\psi}|\mathcal{S}_{n}\right)\right]\\
&&\leq\left[1-P\left(\|\bar{\theta}-\theta_{0}\|^{2}\aplt C_{1}c/card(S_{n})\right)\right]^{card(S_{n})}\\
&&\leq(1-\rho C_{1}/card(S_{n}))^{card(S_{n})}\rightarrow \exp(-\rho C_{1}),
\end{eqnarray*}
where $\rho>0$. The second inequality follows since the cardinality of $\mathcal{S}_{n}$ is larger than $cn^{2\psi}$. The last inequality
follows from the boundedness of the density of $\bar{\theta}$ around $\theta_{0}$. This completes the proof of corollary~\ref{stosea}. $\Box$

{\it Proof of theorem~\ref{thm1}}: We first prove the below lemma~\ref{thm1}.1.

{\it lemma~\ref{thm1}.1.} Assuming the conditions in theorem~\ref{thm1} and that
\begin{eqnarray}
\Pi_{n}(\hat{\theta}_{n}^{(k-1)},t_{n})-\tilde{I}_{0}=O_{P}(r_{n}^{(k-1)}),\label{ass-pi}
\end{eqnarray}
we have
\begin{eqnarray}
(\hat{\theta}_{n}^{(k)}-\hat{\theta}_{n})&=&O_{P}\left(\|\hat{\theta}_{n}^{(k-1)}-\hat{\theta}_{n}\|r_{n}^{(k-1)}\vee \frac{g_{r}(|s_{n}|\vee
n^{-1/2}\vee\|\hat{\theta}_{n}^{(k-1)}-\hat{\theta}_{n}\|)}{n|s_{n}|}\nonumber\right.\\
&&\mbox{\hspace{3in}}\left.\vee|s_{n}| \right)\label{res1}
\end{eqnarray}
for $k=1,2,\ldots$.

 {\it Proof:} Based on (\ref{sche}), we have
\begin{eqnarray}
\Pi_{n}(\hat{\theta}_{n}^{(k-1)},t_{n})\sqrt{n}(\hat{\theta}_{n}^{(k)}-\hat{\theta}_{n})&=&
\left[\sqrt{n}\Pi_{n}(\hat{\theta}_{n}^{(k-1)},t_{n})(\hat{\theta}_{n}^{(k-1)}-\hat{\theta}_{n})\right]+\sqrt{n}\Gamma_{n}
(\hat{\theta}_{n},s_{n})\nonumber\\&+& \left[\sqrt{n}(\Gamma_{n}(\hat{\theta}_{n}^{(k-1)},s_{n})-\Gamma_{n}
(\hat{\theta}_{n},s_{n}))\right].\label{pro1}
\end{eqnarray}
The second term in the above equation equal to
\begin{eqnarray*}
O_{P}\left(\sqrt{n}|s_{n}|\vee \frac{g_{r}(|s_{n}|)}{\sqrt{n}|s_{n}|}\right)
\end{eqnarray*}
according to (\ref{estesco}) and (\ref{lnplexphat}). The third term in (\ref{pro1}) can be written as
\begin{eqnarray*}
-\sqrt{n}\tilde{I}_{0}(\hat{\theta}_{n}^{(k-1)}-\hat{\theta}_{n})+O_{P}\left(\sqrt{n}|s_{n}|\vee
\frac{g_{r}(n^{-1/2}\vee|s_{n}|\vee\|\hat{\theta}_{n}^{(k-1)}-\hat{\theta}_{n}\|)}{\sqrt{n}|s_{n}|}\right).
\end{eqnarray*}
for $k=1,2,\ldots$ by replacing $\hat{\theta}_{n}^{(0)}$ with $\hat{\theta}_{n}$ and $U_{n}$ with $\hat{\theta}_{n}^{(k-1)}-\hat{\theta}_{n}$ in
(\ref{lem2}). Combining the above analysis, the assumption (\ref{ass-pi}) and nonsingularity of $\tilde{I}_{0}$, we complete the proof of
(\ref{res1}). $\Box$

We next start the proof of (\ref{optrate})-(\ref{2-mle}). Combining (\ref{estrel}) with (\ref{res1}), we can obtain that
\begin{eqnarray*}
\hat{\theta}_{n}^{(k)}-\hat{\theta}_{n} &=& O_{P}\left(\frac{g_{r}(|t_{n}|\vee \|\hat{\theta}_{n}^{(k-1)}-
\hat{\theta}_{n}\|)\|\hat{\theta}_{n}^{(k-1)}- \hat{\theta}_{n}\|}{nt_{n}^{2}}\right.\\
&&\mbox{\hspace{1.8in}}\left.\vee|s_{n}|\vee\frac{g_{r}(|s_{n}|\vee
n^{-1/2}\vee\|\hat{\theta}_{n}^{(k-1)}-\hat{\theta}_{n}\|)}{n|s_{n}|}\right)\nonumber\\&=&O_{P}(f_{k-1}(|t_{n}|)\vee g_{k-1}(|s_{n}|)).
\end{eqnarray*}
Considering the form of $g_{r}(\cdot)$ specified in lemma~1, we have $g_{k-1}(|s_{n}|)\geq \tilde{g}(|s_{n}|)\equiv(|s|\vee
n^{-2r-1/2}|s_{n}|^{-1}\vee n^{-3/2}|s_{n}|^{-1})$. The smallest convergence rate for $\tilde{g}(|s_{n}|)$ is $n^{-3/4}(n^{-r-1/4})$ if we
choose $s_{n}\asymp n^{-3/4}$ $(s_{n}\asymp n^{-r-1/4})$. The above analysis implies (\ref{optrate}).

In the below proof of (\ref{1-mle}) and (\ref{2-mle}), we consider different cases when $r\geq 1/2$ and $1/4<r<1/2$, respectively. For $r\geq
1/2$, by some algebra we can show that for $k\geq 1$
\begin{eqnarray}
\hat{\theta}_{n}^{(k)}-\hat{\theta}_{n}=O_{P}(\|\hat{\theta}_{n}^{(k-1)}-\hat{\theta}_{n}\|^{3/2})\label{rll}
\end{eqnarray}
when $\|\hat{\theta}_{n}^{(k-1)}-\hat{\theta}_{n}\|^{-1}=O_{P}(\sqrt{n})$,
$s_{n}^{\ast}\asymp\|\hat{\theta}_{n}^{(k-1)}-\hat{\theta}_{n}\|^{3/2}$ and $t_{n}^{\ast}\asymp\|\hat{\theta}_{n}^{(k-1)}-\hat{\theta}_{n}\|$.
And when $\|\hat{\theta}_{n}^{(k-1)}-\hat{\theta}_{n}\|=O_{P}(n^{-1/2})$, $\|\hat{\theta}_{n}^{(k)}-\hat{\theta}_{n}\|$ achieves the optimal
rate $O_{P}(n^{-3/4})$ given $s_{n}^{\ast}\asymp n^{-3/4}$ and $t_{n}^{\ast}\asymp n^{-1/2}$. Thus we only need to figure out how many iterative
steps needed for $k$-step MLE to achieve root-$n$ rate. From (\ref{rll}), we know that the convergence rate for $N_{1}$-step MLE will be
$O_{P}(n^{-1/2})$, where $N_{1}=int[\log 2\psi/\log(2/3)]$, given $\hat{\theta}_{n}^{(0)}$ is $n^{\psi}$-consistent.  This concludes the proof
for $r\geq 1/2$.

We next show (\ref{2-mle}) when $1/4<r<1/2$. Similarly we have for $k\geq 1$
\begin{eqnarray}
\hat{\theta}_{n}^{(k)}-\hat{\theta}_{n}=O_{P}(\|\hat{\theta}_{n}^{(k-1)}-\hat{\theta}_{n}\|^{3/2})\label{rsl}
\end{eqnarray}
if $\|\hat{\theta}_{n}^{(k-1)}- \hat{\theta}_{n}\|^{-1}=O_{P}(n^{r})$, $s_{n}^{\ast}\asymp\|\hat{\theta}_{n}^{(k-1)}-\hat{\theta}_{n}\|^{3/2}$
and $t_{n}^{\ast}\asymp\|\hat{\theta}_{n}^{(k-1)}-\hat{\theta}_{n}\|$. However if $\|\hat{\theta}_{n}^{(k-1)}-
\hat{\theta}_{n}\|^{-1}=O_{P}(n^{1/2})$ and $\|\hat{\theta}_{n}^{(k-1)}- \hat{\theta}_{n}\|=O_{P}(n^{-r})$, we have for $k\geq 1$
\begin{eqnarray}
\hat{\theta}_{n}^{(k)}-\hat{\theta}_{n}=O_{P}(\|\hat{\theta}_{n}^{(k-1)}-\hat{\theta}_{n}\|^{1/2}n^{-r})\label{rsl0}
\end{eqnarray}
given $s_{n}^{\ast}\asymp\|\hat{\theta}_{n}^{(k-1)}-\hat{\theta}_{n}\|^{1/2}n^{-r}$ and $t_{n}^{\ast}\asymp n^{-r}$. We next consider two-stage
iterations for $K$-step MLE. If $\hat{\theta}_{n}^{(0)}$ is $n^{\psi}$-consistent for $\psi<r$, then at least $M_{1}$ iterations are needed such
that $\|\hat{\theta}_{n}^{(M_{1})}-\hat{\theta}_{n}\|=O_{P}(n^{-r})$ based on (\ref{rsl}), where $M_{1}=int[\log(\psi/r)/\log(2/3)]$. When
$K$-step MLE has achieved the $n^{r}$-consistency, we further need $M_{2}$ steps to achieve the root-$n$ rate, where
$M_{2}=int[\log(4r/(4r-1))/\log(2)-1]$, from (\ref{rsl0}). Then we complete the whole proof for theorem~\ref{thm1}. $\Box$

\vspace{0.2 in} \noindent{\large\bf References}
\begin{description}
\item
Bickel, P. J., Klaassen, C. A. J., Ritov, Y. and Wellner, J. A. (1998). {\it Efficient and Adaptive Estimation for Semiparametric Models}.
Springer-Verlag, New York.
\item
Cheng, G. and Kosorok, M.R. (1987)  (2006). Higher order semiparametric frequentist inference with the profile sampler. {\it Annals of
Statistics}, Accepted.

\item
Cheng, G. and Kosorok, M.R. (2007). General Frequentist Properties of the Posterior Profile Distribution. {\it Annals of Statistics}, Invited
Revision.

(http://arxiv.org/abs/math.ST/0612191)

\item
Cox, D. R. (1972). Regression models and life-tables. {\em Journal of the Royal Statistical Society, Series B} {\bf 34} 187--220.

\item
Fan, J. and Chen, J. (1999). One-Step Local Quasi-Likelihood Estimation. {\it Journal of the Royal Statistical Society, Series B} {\bf 61}
927--943.

\item
Fan, J., Lin, H. and Zhou, Y. (2006). Local Partial-Likelihood Estimation for Lifetime Data. {\it Annals of Statistics} {\bf 34} 290--325.

\item
 Groeneboom, P. (1991). Nonparametric maximum likelihood estimators for interval censoring and deconvolution. {\it Technical Report} {\bf 378},
Department of Statistics, Stanford University.

\item
Huang, J. (1996). Efficient estimation for the Cox model with interval censoring. {\it Annals of Statistics} {\bf 24,} 540--568.

\item
Jassen, P., Jureckova, J. and Veraverbeke, N. (1985). Rate of Convergence of One- and Two-stpe M-estimators with Applications to Maximum
Likelihood and Pitman Esitmators. {\it Annals of Statistics} {\bf 25} 1471--1509.

\item
Kosorok, M. R., Lee, B. L. and Fine, J. P. (2004). Robust inference for univariate proportional hazards frailty regression models. {\it Annals
of Statistics} {\bf 32} 1448--1491.

\item
 Kosorok, M. R. (2007). {\em Introduction to Empirical Processes and Semiparametric Inference}. Springer, New York.

\item
 Lee, B. L., Kosorok, M. R. and Fine, J. P. (2005). The profile sampler. {\it Journal of the American Statistical Association} {\bf 100}
960--969.

\item
Murphy, S. A. and Van der Vaart, A. W. (1999). Observed information in semiparametric models. {\em Bernoulli} {\bf 5} 381--412.

\item
Murphy, S. A. and Van der Vaart, A. W. (2000). On profile likelihood. {\it Journal of the American Statistical Association} {\bf 93} 1461--1474.

\item
van der Vaart, A. W. (1998). {\em Asymptotic Statistics}. Cambridge University Press, Cambridge.

\item
Cai, Z., Fan, J. and Li, R. (2000). Efficient Estimation and Inferences for Varying-Coefficient Models. {\it Journal of the American Statistical
Association} {\bf 95} 888--902.
\end{description}

\vskip .65cm \noindent Department of Statistical Science\\ Duke University\\ Durham, NC, U.S. \vskip 2pt \noindent E-mail: chengg@stat.duke.edu
\\\noindent Phone: 919-684-5956 \\\noindent Fax: 919-684-8594

\vskip 2pt \newpage

\begin{center}
Table 1. {\it Cox regression with right censored data($\theta_{0}=1$ and $500$ samples).} 
\centering 
\begin{tabular}{c c c c c c} 
\hline\hline 
n & $\hat{\theta}_{n}^{(0)}$ & $\hat{\theta}_{n}^{(1)}$ & $\hat{\theta}_{n}$ & $n^{3/4}|\hat{\theta}_{n}-\hat{\theta}_{n}^{(1)}|$\\
\hline 
50 & 1.0229 & 1.0222 & 1.0167 & 0.1030 \\
100 & 1.0346 & 1.0344 & 1.0324 & 0.0632 \\
200 & 0.9979 & 0.9979 & 1.0028 & 0.2606 \\
500 & 0.9974 & 0.9974 & 0.9964 & 0.1057 \\
\multicolumn{6}{@{}p{12.6cm}@{}}{\rule{12.6cm}{0.2pt}}\\
\end{tabular}
\end{center}

\begin{center}
Table 2. {\it  Cox regression with current status data ($\theta_{0}=1$ and $500$ samples).} 
\centering 
\begin{tabular}{c c c c c c c} 
\hline\hline 
n & $\hat{\theta}_{n}^{(0)}$ & $\hat{\theta}_{n}^{(1)}$ & $\hat{\theta}_{n}^{(2)}$ & $\hat{\theta}_{n}^{(3)}$ & $\hat{\theta}_{n}$ &
$n^{7/12}|\hat{\theta}_{n}-\hat{\theta}_{n}^{(3)}|$\\
[0.5ex] 
\hline 
50 & 1.0452 & 1.8218 & 1.7226 & 1.7563 & 1.1962 & 5.4870\\
100 & 0.8017 & 0.7604 & 0.7997 & 0.8289 & 0.8541 & 0.3699\\
200 & 0.8118 & 0.7692 & 0.8474 & 0.8425 & 0.8859 & 0.9545\\
500 & 0.8376 & 0.8364 & 0.8757 & 0.9592 & 0.9896 & 1.1410 \\
\multicolumn{7}{@{}p{12.6cm}@{}}{\rule{12.6cm}{0.2pt}}\\
\multicolumn{7}{@{}p{12.6cm}@{}}
{\small%
$n$, sample size; $\hat{\theta}_{n}^{(0)}$, initial estimate; $\hat{\theta}_{n}^{(1)}$, one step MLE; $\hat{\theta}_{n}^{(2)}$, two step MLE;
$\hat{\theta}_{n}^{(3)}$, three step MLE; $\hat{\theta}_{n}$, MLE.}
\end{tabular}
\end{center}

\vskip .3cm


\end{document}